\documentclass[a4paper,11pt]{article}
\usepackage[margin=2.5cm]{geometry}

\usepackage{mathptmx}       
\usepackage{helvet}         
\usepackage{courier}        
\usepackage{type1cm}        
\usepackage{makeidx}         
\usepackage{graphicx}        
\usepackage{multicol}        
\usepackage[bottom]{footmisc}
\usepackage{amsmath, amsfonts, amssymb}
\usepackage{color}
\usepackage{tikz}
\usepackage{todonotes}
\usepackage{bm}
\usepackage{hyperref} 

\newcommand*\samethanks[1][\value{footnote}]{\footnotemark[#1]}

\newtheorem{theorem}{Theorem}

\newtheorem{remark}[theorem]{Remark}
\numberwithin{equation}{section}
\numberwithin{theorem}{section}

\newcommand{\divergence}{\boldsymbol \nabla \cdot}
\newcommand{\gradient}{\boldsymbol \nabla}
\newcommand{\llangle}{\left\langle}
\newcommand{\rrangle}{\right\rangle}
\DeclareMathOperator*{\argmin}{\mathrm{arg\,min}}

\newcommand{\strain}{\boldsymbol \varepsilon}
\newcommand{\displacement}{\boldsymbol u}
\newcommand{\flux}{\boldsymbol q}
\newcommand{\pressure}{p}
\newcommand{\fluidvolume}{\theta}

\newcommand{\permeability}{\boldsymbol \kappa}
\newcommand{\viscosity}{\mu}

\newcommand{\f}{\boldsymbol f}

\newcommand{\displacementSpace}{\bm{\mathcal{V}}}
\newcommand{\pressureSpace}{\mathcal{Q}}
\newcommand{\fluxSpace}{\bm{\mathcal{W}}}

\newcommand{\testdisplacement}{\boldsymbol v}
\newcommand{\testpressure}{s}
\newcommand{\testflux}{\boldsymbol w}

\makeindex             

\begin{document}

\title{Free energy diminishing discretization of Darcy-Forchheimer flow in poroelastic media}

\author{Jakub W.\ Both \thanks{Department of Mathematics, University of Bergen, Norway} \ \thanks{Corresponding author: jakub.both@uib.no}
\and Jan M.\ Nordbotten\samethanks[1]
\and Florin A.\ Radu\samethanks[1]
} 

\date{}


\maketitle

\abstract{In this paper, we develop a discretization for the non-linear coupled model of classical Darcy-Forchheimer flow in deformable porous media, an extension of the quasi-static Biot equations. The continuous model exhibits a generalized gradient flow structure, identifying the dissipative character of the physical system. The considered mixed finite element discretization is compatible with this structure, which gives access to a simple proof for the existence, uniqueness, and stability of discrete approximations. Moreover, still within the framework, the discretization allows for the development of finite volume type discretizations by lumping or numerical quadrature, reducing the computational cost of the numerical solution.
%
}


\section{Introduction}

Flow in deformable porous media has been of increased interest in the recent past. Applications of societal and industrial relevance range from geotechnical to biomedical engineering, including the consolidation of the subsurface due to fluid production and the deformation of fluid-filled soft tissue.

Regarding slow viscous flow in linearly poroelastic media, the quasi-static linear Biot equations are often chosen as mathematical model, essentially coupling equations of linear elasticity and single-phase flow. For applications with significantly faster flow rates, Darcy's law is not further applicable. Instead the classical non-linear Darcy-Forchheimer law~\cite{Forchheimer1901} is often utilized, cf., e.g.,~\cite{Markert2007,Atik2019}.

In this paper, we study the discretization of Darcy-Forchheimer flow in poroelastic media. The basis for this will be a (mixed) finite element method -- widely used in the context of flow in porous media since being locally mass conservative. However, typically it suffers from larger algebraic systems to be solved, compared to, e.g., cell-centered finite volume discretizations. To circumvent this, for Darcy flow in non-deformable media, mass lumping~\cite{Baranger1996} or approximate numerical quadrature techniques resulting in the (symmetric) multipoint flux mixed finite element method~\cite{Wheeler2006} have been developed allowing for local discrete flux pressure relationships. These are related to finite volume schemes employing respectively two-point and multipoint flux approximations~\cite{Baranger1996,Klausen2006,Klausen2008}. Moreover, the resulting linear systems involve block-diagonal mass matrices, which allow for efficient solution. Recently, these techniques have been also applied in the context of deformable media~\cite{Hu2017,Ambartsumyan2020}. Regarding Darcy-Forchheimer flow in non-deformable media, especially the mixed finite element method on unstructured meshes~\cite{Park2005,Girault2008} and (similar to the above efforts) a block-centered finite difference method on rectangular grids~\cite{Rui2012} have been developed and studied in more detail including their well-posedness and theoretical convergence.  

Motivated by all these advances, we propose a combination of the mixed finite element method and similar localization techniques in the context of Darcy-Forchheimer flow in deformable media. We particularly emphasize the inherent gradient flow structure of the continuous model, quantifying the dissipation of free energy over time. By construction, the numerical schemes considered here mimic a similar structure. Remarkably, it gives access to simple well-posedness and stability analyses of the numerical schemes.

The outline of the remaining paper is as follows. In section~\ref{section:model}, the mathematical model is described. In section~\ref{section:discretization}, the numerical method is presented, for which theoretical properties are discussed in section~\ref{section:analysis}. 

Not part of this paper, but in the future, a numerical study will be conducted with focus on assessing the potential accuracy loss of the localization techniques, and efficiency gain regarding the algebraic solution. In particular, the exploitation of the block-diagonal nature of the flux mass term will be combined with robust splitting schemes as in~\cite{Both2017,Both2019}, benefiting from the linear character of the elasticity equations.

\section{Model for Darcy-Forchheimer flow in poroelastic media}
\label{section:model}

The mathematical model couples the balance of linear momentum and the conservation of mass for a poroelastic medium, here modeled as an open, connected domain $\Omega\subset\mathbb{R}^d$, $d\in\{2,3\}$. In addition, constitutive relations are considered: the medium is assumed to be linearly elastic and fully saturated with a slightly compressible fluid, with fluid flow described by the classical Darcy-Forchheimer law~\cite{Forchheimer1901}. The solid-fluid interaction is governed by the so-called effective stress. Finally, the system of governing equations reads \\[-0.75cm]

\begin{subequations}
\label{model:continuous}
\begin{alignat}{2}
  -\divergence \left( \mathbb{C}\strain(\displacement) - \alpha \pressure \,\mathbf{I} \right) &= \f &&\text{on }\Omega,\\
 \partial_t \left( \frac{1}{M} p + \alpha \divergence \displacement \right) + \divergence \flux &= 0 &&\text{on }\Omega,\\
 \viscosity \permeability^{-1} \flux + \rho \beta |\flux| \flux &= -\gradient p &\hspace{0.5cm}&\text{on }\Omega,
 \label{model:continuous:darcy}
\end{alignat}
\end{subequations}
where the primal variables are the displacement $\displacement:\Omega \rightarrow \mathbb{R}^d$, the fluid pressure $p:\Omega \rightarrow \mathbb{R}$, and the volumetric flux $\flux: \Omega \rightarrow \mathbb{R}^d$. Moreover, $\mathbb{C}$ denotes the (symmetric positive definite) fourth-order stiffness tensor, $\strain(\cdot)$ denotes the linear strain tensor, $\alpha\in(0,1]$ is the Biot coefficient, $\f$ is a an external force, $\partial_t$ denotes the derivative in time, $M\geq 0$ is the modulus accounting for the compressibility of the fluid and solid grains, $\viscosity>0$ is the fluid viscosity, $\permeability$ is the (symmetric positive definite) permeability tensor, $\rho>0$ is the reference fluid density, and $\beta\geq 0$ is the Forchheimer index; the case $\beta=0$ simplifies to Darcy's law. Ultimately,~\eqref{model:continuous:darcy} can be viewed as non-linear Darcy law with a direction dependent mobility. 

For the sake of brevity, all material parameters are assumed to be constant. Furthermore, no external body or surface sources for the volume content are considered. However, we note that corresponding extensions are possible.

The system is closed with initial conditions $\fluidvolume(0)=\fluidvolume^0$ for the volume content $\fluidvolume=\tfrac{1}{M}p + \alpha \divergence \displacement$ as well as boundary conditions for the displacement and flux. Here, for simplicity we choose ${\displacement_|}_{\partial \Omega} = \bm{0}$ and ${(\flux \cdot \bm{n})_|}_{\partial \Omega} = 0$, where $\bm{n}$ denotes the outward normal on $\partial\Omega$.

\subsection{The gradient flow structure of the model} 

As presented in~\cite{Both2019}, the weak formulation of model~\eqref{model:continuous} exhibits a generalized gradient flow structure, in the sense of the lecture notes by Peletier~\cite{Peletier2014}. Short, one can define the standard poroelastic Helmholtz free energy and a non-quadratic dissipation potential extending the classical potential corresponding to linear Darcy flow
\begin{align*}
 \mathcal{E}(\displacement,\fluidvolume) 
 &= \frac{1}{2} \llangle \mathbb{C} \strain(\displacement), \strain(\displacement) \rrangle 
 + \frac{M}{2} \left\| \theta - \alpha \divergence \displacement \right\|_{L^2(\Omega)}^2
 -
 \llangle \f, \displacement \rrangle, \\
 \mathcal{D}(\flux)
 &= \frac{\viscosity}{2} \llangle \permeability^{-1} \flux, \flux \rrangle + \frac{\beta \rho}{3} \left\| \flux \right\|_{L^3(\Omega)}^3,
\end{align*}
with $\|\cdot\|_{L^p(\Omega)}$ and $\llangle \cdot, \cdot \rrangle$ denoting the standard $L^p(\Omega)$ norm and $L^2(\Omega)$ scalar product, respectively. In the following, $H^1(\Omega)^d$ denotes the Sobolev space consisting of $L^2(\Omega)^d$ functions with weak derivatives in $L^2(\Omega)^{d \times d}$, and $H(\mathrm{div};\Omega)$ requires solely the divergence to be in $L^2(\Omega)$. Moreover, let $\frac{\delta}{\delta (\cdot)}$ denote the Fr\'echet differential.

In the absence of dissipation of energy due to solid deformation, weak solutions to~\eqref{model:continuous} are alternatively characterized by the degenerate generalized gradient flow
\begin{align*}
 \displacement &= \argmin_{\testdisplacement \in H^1(\Omega)^d} \, \mathcal{E}(\testdisplacement,\fluidvolume),\\
 (\partial_t\fluidvolume, \flux) &= \argmin_{(\testpressure,\testflux) \in L^2(\Omega) \times L^3(\Omega)^d \cap H(\mathrm{div};\Omega)}\, \left\{ \mathcal{D}(\testflux) + \llangle \frac{\delta \mathcal{E}}{\delta \theta} (\displacement,\fluidvolume), \testpressure \rrangle\right\}\\[0.2cm]
 &\qquad \text{subj.\ to}  \qquad
\left\{ \begin{array}{rll}
\testdisplacement &= \bm{0} &\text{ on }\partial\Omega,\\
 \testpressure + \divergence \testflux &= 0, &\text{ on }\Omega,\\
 \testflux \cdot \bm{n} & = 0 &\text{ on }\partial \Omega,
 \end{array}\right. 
\end{align*}
such that the flux $\flux$ governs $\partial_t \fluidvolume$. 
The governing equations~\eqref{model:continuous} can be recovered as the optimality conditions.
The fluid pressure $p$ enters as Lagrange variable associated to mass conservation, as well as dual variable $p = \frac{\delta \mathcal{E}}{\delta \fluidvolume}
(\displacement,\fluidvolume)$ and $p = \frac{\delta \mathcal{E}}{\delta \divergence \displacement}(\displacement,\fluidvolume)$.

For sufficiently smooth solutions, employing the chain rule and the convexity of $\mathcal{D}$ yield the following energy--dissipation relation
\begin{align*}
 \frac{d}{dt} \mathcal{E} (\displacement,\theta)
 =
 \llangle \frac{\delta \mathcal{E}(\displacement, \fluidvolume)}{\delta (\displacement,\fluidvolume)}, (\partial_t \displacement,\partial_t \fluidvolume) \rrangle
 =
 -\llangle p, \divergence \flux \rrangle
 =
 - \llangle \frac{\delta \mathcal{D}}{\delta \flux}(\flux), \flux \rrangle
 \leq 
 - \mathcal{D}(\flux).
 \end{align*} 

\begin{remark}[Incompressible case\label{remark:incompressible-case-model}]
 In the incompressible case, i.e., $M=\infty$, the energy degenerates and becomes merely sub-differentiable, as then (here for $\f=\bm{0}$)
 \begin{align*}
  \mathcal{E}(\displacement,\fluidvolume) = \frac{1}{2} \llangle \mathbb{C} \strain(\displacement), \strain(\displacement) \rrangle + \chi\left(\|\fluidvolume - \alpha \divergence \displacement \|_{L^2(\Omega)}\right),\quad\text{with}\quad \chi(p)=\left\{\begin{array}{ll} 0, &p=0,\\ \infty, &\text{else.} \end{array} \right.
 \end{align*}

\end{remark}

\section{Numerical discretization\label{section:discretization}}

We present discretizations mimicking the free energy dissipating character of the continuous problem formulation. Focussing on the non-linear character of the problem, we limit the discussion to equidistant time-stepping and simplicial grids. 

\subsection{Semi-discrete approximation in variational form}

For the discretization in time, we utilize the minimizing movement scheme, which for~\eqref{model:continuous} is equivalent with the implicit Euler method. Given a time interval of interest $[0,T]$, with $T>0$ denoting the final time, we consider an equidistant partition $0=t_0 < t_1 < \ldots < t_N =T$ with time step size $\tau$. Furthermore, let $\mathcal{J}(\displacement, \fluidvolume, \flux):=  \tau \mathcal{D}(\flux) + \mathcal{E}(\displacement, \fluidvolume)$. The discretization at time step $n\geq 1$ reads: 
\textit{given $\fluidvolume^{n-1}\in L^2(\Omega)$, find $(\displacement^n,\fluidvolume^n,\flux^n)\in H^1(\Omega)^d \times L^2(\Omega) \times L^3(\Omega)^d \cap H(\mathrm{div};\Omega)$ satisfying}
\begin{align*}
 (\displacement^n,\fluidvolume^n,\flux^n) = \argmin_{(\displacement,\fluidvolume,\flux)}\, \mathcal{J}(\displacement, \fluidvolume, \flux),
 \ 
 \text{subj.\ to}
 \ 
 \left\{
 \begin{array}{rll}
 \displacement &= \bm{0}, &\text{ on }\partial\Omega,\\[-0.1cm]
 \frac{\fluidvolume - \fluidvolume^{n-1}}{\tau} + \divergence \flux &= 0, & \text{ on }\Omega,\\
 \flux \cdot \bm{n} &= 0, &\text{ on }\partial\Omega.
 \end{array}\right. 
\end{align*}

\subsection{Fully discrete approximation in variational form}

For the discretization in space, we utilize the (mixed) finite element method with the possibility for reduced computational complexity using lumping or appropriate numerical quadrature. Starting with the semi-discrete formulation, we introduce a mesh-dependent version $\mathcal{J}_h$ of the original objective function $\mathcal{J}$ for this. \\[-0.25cm]

\noindent
\textbf{The mesh.} Assume the physical domain $\Omega$ is polygonal and can be partitioned by simplices. Let $\mathcal{T}_h$ denote such a simplicial mesh with elements $K\in \mathcal{T}_h$ and faces $e\in \mathcal{F}_h$; 
let $\{\bm{r}_i\}_{i=1,..,d+1}$ denote the corners of $K\in\mathcal{T}_h$.\\[-0.25cm] 

\noindent
\textbf{The finite element spaces.} We consider classical conforming approximation spaces (including essential boundary conditions): $\displacementSpace_h$, piecewise linear elements for the displacements; $\pressureSpace_h$, piecewise constant elements for the fluid content and the fluid pressure; and either $\fluxSpace_h^\mathtt{RT}$, lowest order Raviart-Thomas, or $\fluxSpace_h^\mathtt{BDM}$, lowest order Brezzi-Douglas-Marini elements, for the flux, depending on the subsequent choice for $\mathcal{J}_h$. If the particular choice is not crucial, we write $\fluxSpace_h$ and allow for $\fluxSpace_h^\mathtt{RT}$ and $\fluxSpace_h^\mathtt{BDM}$. For detailed introduction of the finite element spaces, we refer to~\cite{Boffi2013}.\\[-0.25cm]

\noindent
\textbf{Objective function for MFEM.} Instead of pursuing the Galerkin method for deriving fully discrete approximations, a mesh-dependent approximation of the objective function $\mathcal{J}$ is additionally considered
\begin{align*}
  \mathcal{J}_h(\displacement_h,\fluidvolume_h,\flux_h):= \tau \mathcal{D}_h(\flux_h)
  + \mathcal{E}_h(\displacement_h,\fluidvolume_h).
\end{align*}
The canonical mixed finite element discretization of~\eqref{model:continuous} results in particular for
\begin{align*}
 \mathcal{D}_h(\flux_h) &:= \mathcal{D}(\flux_h),\\
 \mathcal{E}_h(\displacement_h, \fluidvolume_h) &:= \frac{1}{2} \llangle \mathbb{C} \strain(\displacement_h), \strain(\displacement_h) \rrangle + \frac{M}{2} \left\| \Pi_{\pressureSpace_h} \left( \fluidvolume_h  - \alpha \divergence \displacement_h \right) \right\|_{L^2(\Omega)}^2,
\end{align*}
for $\flux_h\in \fluxSpace_h$ and $(\displacement_h,\fluidvolume_h)\in \displacementSpace_h \times \pressureSpace_h$, where $\Pi_{\pressureSpace_h}$ denotes the $L^2(\Omega)$ projection onto $\pressureSpace_h$, allowing for measuring the fluid energy in the 'units' of the fluid volume. In the incompressible case ($M=\infty$), $\mathcal{E}_h$ is defined using an indicator function, as $\mathcal{E}$. \\[-0.25cm]

\noindent
\textbf{Definition of the method.} Given a suitable approximation $\fluidvolume_h^0\in \pressureSpace_h$ of the initial datum $\fluidvolume^0 \in L^2(\Omega)$, the fully discrete approximation at time step $n\geq 1$ reads: \textit{given $\fluidvolume_h^{n-1}\in \pressureSpace_h$, find $(\displacement_h^n,\fluidvolume_h^n,\flux_h^n)\in \displacementSpace_h \times \pressureSpace_h \times \fluxSpace_h$ satisfying}\\[-0.75cm]

\begin{subequations}
\label{discretization-minimization}
\begin{align}
\label{discretization-minimization-1}
 &(\displacement_h^n,\fluidvolume_h^n,\flux_h^n) 
 = \argmin_{(\displacement_h,\fluidvolume_h,\flux_h)}\, \tau \mathcal{D}_h(\flux_h) + \mathcal{E}_h(\displacement_h,\fluidvolume_h)\\
\label{discretization-minimization-2}
 &\qquad\text{subj.\ to}
 \ \ 
 |K|\left(\fluidvolume_K - \fluidvolume_K^{n-1} \right) + \tau\int_{K} \divergence \flux_h \, dx = 0 \quad \forall K\in\mathcal{T}_h.
\end{align}
\end{subequations}

\noindent
\textbf{Saddle point formulation.} The optimality conditions corresponding to~\eqref{discretization-minimization} are obtained by introducing a Lagrange multiplier, which eventually can be identified as the discrete fluid pressure $p_h\in \pressureSpace_h$.  We skip the calculations and directly present the final discrete system to be solved at time step $n\geq 1$: \textit{given $(\displacement_h^{n-1},\pressure_h^{n-1})\in \displacementSpace_h \times \pressureSpace_h$, find $(\displacement_h^n,p_h^n,\flux_h^n)\in \displacementSpace_h \times \pressureSpace_h \times \fluxSpace_h$ satisfying for all test functions $(\testdisplacement_h,\testflux_h)\in \displacementSpace_h \times \fluxSpace_h$ and elements $K\in \mathcal{T}_h$}
\begin{subequations} \label{discretization-saddle-point}
\begin{align}
\llangle \mathbb{C} \strain(\displacement_h^n), \strain(\testdisplacement_h) \rrangle - \alpha \llangle p_h^n, \divergence \testdisplacement_h \rrangle &= \llangle \f, \testdisplacement_h \rrangle,\\[0.3cm]
\label{discretization-darcy}
m(\flux_h^n;\flux_h^n,\testflux_h) - \llangle p_h^n, \divergence \testflux_h \rrangle &= 0,\\[0.1cm]
\label{discretization-pressure}
\frac{|K|}{M} (\pressure_K^n - \pressure_K^{n-1}) + \alpha \int_K \divergence ( \displacement_h^n - \displacement_h^{n-1} )\,dx + \tau \int_{K} \divergence \flux_h^n \, dx &= 0,
\end{align}
\end{subequations}
where the non-linear form is given by $m(\bm{u};\bm{v},\bm{w})
 =
\int_\Omega \left(\viscosity \permeability^{-1} + \rho \beta |\bm{u}| \mathbf{I} \right) \bm{v} \cdot \bm{w} \, dx$.
The volume content at time step $n$ can be recovered by $\fluidvolume_h^n = \tfrac{1}{M} \pressure_h^n + \alpha \Pi_{\pressureSpace_h} \divergence \displacement_h^n$.\\[-0.25cm]

\noindent
\textbf{Localized dissipation potential.} One drawback of the above formulation is the non-local relation between fluxes and pressure gradients. From a computational point of view, this results in an involved numerical solution. Instead, motivated by efforts for linear problems, we consider two choices: (i) mass lumping as in~\cite{Baranger1996,Hu2017} resulting in a two-point flux type approximation for $\flux_h \in \fluxSpace_h^\mathtt{RT}$
\begin{align*}
\mathcal{D}_h^\mathtt{ML}(\flux_h) &= \sum_{e\in \mathcal{F}_h} \omega_e\left[
 \frac{\mu}{2\kappa} \left| \int_e \flux_h \cdot \bm{n}_e \, ds\right|^2 +
 \frac{\beta}{3}  \left| \int_e \flux_h \cdot \bm{n}_e \, ds\right|^3 \right]
\end{align*}
where $\bm{n}_e$ is a uniquely defined normal on $e\in\mathcal{F}_h$ and $\omega_e$ is a suitable weight involving distances and measures of geometric entities~\cite{Baranger1996}; and (ii) trapezoidal quadrature as in~\cite{Wheeler2006} resulting in a multipoint flux type approximation for $\flux_h \in \fluxSpace_h^\mathtt{BDM}$
\begin{align*}
\mathcal{D}_h^\mathtt{Q}(\flux_h) &= \sum_{K \in \mathcal{T}_h}
 \frac{|K|}{d+1} \sum_{i=1}^{d+1} \left[ \frac{\viscosity}{2} {\permeability_|}_K^{-1}(\bm{r}_i) {{\flux_h}_|}_K(\bm{r}_i) \cdot {{\flux_h}_|}_K(\bm{r}_i) + \frac{\beta\rho}{3} \big|  {{\flux_h}_|}_K(\bm{r}_i) \big|^3 \right].
\end{align*}
Lumping is only a sufficient approximation for scalar permeabilities~\cite{Baranger1996}. Also one has to be aware of the fact that only the normal component of the flux contributes here, which is inconsistent with the constitutive Darcy-Forchheimer law.

The corresponding optimality conditions read as~\eqref{discretization-saddle-point} but with an approximation $m_h$ of $m$, which after suitable linearization results in a block-diagonal matrix. Thereby, fluxes may be explicitly represented in terms of pressure values at cell centers. This finally enables the development of efficient numerical solvers.

\section{Existence, uniqueness, and stability}
\label{section:analysis}

The scheme~\eqref{discretization-saddle-point} satisfies local mass conservation at nonlinear and linearized level -- independent of the particular choices for $\fluxSpace_h$ and $\mathcal{D}_h$. In the following, we additionally emphasize a free energy dissipating character, naturally deduced the inherited minimization structure~\eqref{discretization-minimization}. This structure in particular also guarantees unique solutions to~\eqref{discretization-saddle-point}. For the main theoretical result, we require consistency of the original and approximate dissipation potentials; motivated by~\cite{Baranger1996,Wheeler2006}, we expect that there exists a broad class of grids and parameters
for which the following holds.\\[-0.25cm]

\noindent
\textbf{Assumption 1.}  There exist constants $0<c\leq C <\infty$ such that  $c \,  \mathcal{D}(\flux_h) \leq \mathcal{D}_h(\flux_h) \leq C \,\mathcal{D}(\flux_h)$ for all $\flux_h\in \fluxSpace_h$.

\begin{theorem}[Existence, uniqueness and stability]
 Let $\mathcal{D}_h$ satisfy Assumption 1. Let $M<\infty$ and the remaining material parameters be as introduced in Section~\ref{section:model}. In addition, let $m_h^0\in \pressureSpace_h$ and $\displacement_h^0\in \displacementSpace_h$, then for all time steps $n\!\geq\! 1$, the numerical scheme~\eqref{discretization-saddle-point} has a unique solution $(\displacement_h^n,\pressure_h^n,\flux_h^n)\in\displacementSpace_h \times \pressureSpace_h \times \fluxSpace_h$. It furthermore satisfies the following discrete energy dissipation inequality: for all $N>0$ it holds
 \begin{align*}
  \llangle \mathbb{C}\strain(\displacement_h^N), \strain(\displacement_h^N) \rrangle 
  +
  \frac{1}{M} \left\| p_h^N \right\|^2 + \sum_{n=0}^N \tau \mathcal{D}_h(\flux_h^n) \leq \llangle \mathbb{C}\strain(\displacement_h^{0}), \strain(\displacement_h^{0}) \rrangle 
  +
  \frac{1}{M} \left\| p_h^{0} \right\|^2.
 \end{align*}
\end{theorem}

The proof follows by induction with the induction step for $n\geq 1$ reading as follows.
The existence and uniqueness result follows from classical convex analysis~\cite{Ekeland1999} applied to the minimization formulation~\eqref{discretization-minimization}. 

The discrete function spaces $\displacementSpace_h$, $\pressureSpace_h$, $\fluxSpace_h$ are reflexive Banach spaces, equipped with natural norms: the $H^1(\Omega)$ semi-norm $\|\cdot\|_{\displacementSpace}:=|\cdot|_{H^1(\Omega)}$ on $\displacementSpace_h$, the $L^2(\Omega)$ norm $\|\cdot\|_{\mathcal{Q}}:=\|\cdot\|_{L^2(\Omega)}$ for $\pressureSpace_h$, and $\|\cdot\|_{\fluxSpace}$ on $\fluxSpace_h$, defined by
 \begin{align*}
  \|\testflux_h\|_{\fluxSpace}:= \|\testflux_h\|_{L^2(\Omega)} + \beta\|\testflux_h\|_{L^3(\Omega)} + \| \divergence \testflux_h \|_{L^2(\Omega)},\qquad \testflux_h\in\fluxSpace_h.
 \end{align*}
 The objective function $\mathcal{J}_h$ is strictly convex on $\displacementSpace_h \times \pressureSpace_h \times \fluxSpace_h$, which follows directly from the definition and Assumption 1 ensuring the definiteness of $\mathcal{D}_h$; we emphasize that $\|\fluidvolume_h\|^2$ can be isolated by hiding the coupling term in the elastic energy. Furthermore, $\mathcal{J}_h$ is lower semi-continuous, and proper on
 \begin{align*}
  C:=\left\{ (\displacement_h,\fluidvolume_h,\flux_h) \in \displacementSpace_h \times \pressureSpace_h \times \fluxSpace_h \, \big| \, \text{ \eqref{discretization-minimization-2} is satisfied } \right\}
 \end{align*}
since $\mathcal{J}_h(\bm{0},\fluidvolume_h^{n-1},\bm{0})<\infty$ with $(\bm{0},\fluidvolume_h^{n-1},\bm{0})\in C$. The feasible set $C$ is by that not only convex and closed, but also non-empty. Lastly, again under Assumption 1, $\mathcal{J}_h$ is coercive over $C$; for this, we particularly note that $(\displacement_h,\pressure_h,\flux_h)\in C$ with $\|\divergence \flux_h\|_{L^2(\Omega)} \rightarrow \infty$ implies $\| \fluidvolume_h \|_{L^2(\Omega)} \rightarrow \infty$, which eventually results in $J(\displacement_h,\fluidvolume_h,\flux_h)\rightarrow \infty$. Ultimately, existence and uniqueness of a solution $(\displacement_h^n,\fluidvolume_h^n,\flux_h^n)\in C$ to~\eqref{discretization-minimization} follows from classical convex analysis~\cite{Ekeland1999}. By comparing the solution with the feasible competitor $(\displacement_h^{n-1},\fluidvolume_h^{n-1},\bm{0}) \in C$, stability can be concluded
\begin{align*}
  \mathcal{J}_h(\displacement_h^n,\fluidvolume_h^n,\flux_h^n) = \mathcal{E}_h(\displacement_h^n,\fluidvolume_h^n) + \tau \mathcal{D}_h(\flux_h^n) \leq \mathcal{E}_h(\displacement_h^{n-1},\fluidvolume_h^{n-1}) = \mathcal{J}_h(\bm{0},\fluidvolume_h^{n-1},\bm{0}).
 \end{align*}
 
 Existence, uniqueness, and stability of a solution $(\displacement_h^n,p_h^n,\flux_h^n) \in \displacementSpace_h \times \pressureSpace_h \times \fluxSpace_h$ to the saddle point formulation~\eqref{discretization-saddle-point} then follows from the equivalence of~\eqref{discretization-minimization} and~\eqref{discretization-saddle-point} and the identification $\fluidvolume_h^n = \tfrac{1}{M} p_h^n + \alpha \Pi_{\pressureSpace_h} \divergence \displacement_h^n$. The equivalence follows due to the well-known inf-sup stability of $\fluxSpace_h \times \pressureSpace_h$~\cite{Boffi2013}.

\begin{remark}[Incompressible case]
 In the case of an incompressible medium, i.e., $M=\infty$, both problem formulations~\eqref{discretization-minimization} (modified by the indicator function, cf.\ Rem.~\ref{remark:incompressible-case-model}) and~\eqref{discretization-saddle-point} are still equivalent, if also $\displacementSpace_h \times \pressureSpace_h$ is inf-sup stable wrt.\ the divergence operator. This is not the case for $\displacementSpace_h$ as defined above~\cite{Boffi2013}.
\end{remark}

\begin{remark}[Convergence] 
The convergence of the numerical approximation towards the continuous solution for decreasing mesh and time step size is a delicate subject - in particular regarding the localization techniques. For instance, mass lumping (similar to linear TPFA) is well-known to be prone for errors even in the linear case, e.g., for anisotropic permeability. A further study will be conducted in the future.
\end{remark}

\subsubsection*{Acknowledgement}

This work is supported in part by the Research Council of Norway Project 250223, as well as the FracFlow project funded by Equinor through Akademiaavtalen. 
\vspace{-0.25cm}


\end{document}